\documentclass[12pt,twoside]{paper}
\usepackage{graphicx,fancyhdr,url,amsfonts}
\usepackage[english]{babel}
\usepackage{algorithm}
\usepackage{algorithmic}

\newtheorem{teo}{Theorem}
\newtheorem{prop}{Proposition}
\newtheorem{defi}{Definition}
\newtheorem{obs}{Remark}

\begin{document}
\fancyhead[RO]{H. Oviedo}

 \author{Harry F. Oviedo
        \thanks{Departamento de Ciencias de Computaci\'on,
     Centro de Investigaci\'on en Matem\'aticas, CIMAT A.C. Guanajuato, Guanajuato, Mexico.}}
 \title{A Spectral Gradient Projection Method for the Positive Semi-definite Procrustes Problem}
 \shortauthor{H. Oviedo}

\maketitle

 \begin{abstract}
This paper addresses the positive semi-definite procrustes problem (PSDP). The PSDP corresponds to a least squares problem over the set of
symmetric and semi-definite positive matrices. These kinds of problems appear in many applications such as structure analysis, signal processing, among others. A
non-monotone spectral projected gradient algorithm is proposed to obtain a numerical solution for the PSDP. The proposed algorithm employs the Zhang and Hager's
non-monotone technique in combination with the Barzilai and Borwein's step size to accelerate convergence. Some theoretical results are presented. Finally,
numerical experiments are performed to demonstrate the effectiveness and efficiency of the proposed method, and comparisons are made with other state-of-the-art
algorithms.
\end{abstract}


\section{Introduction}
\label{sec:intro}
The positive semi-definite Procrustes problem (PSDP) is defined as follows: given two rectangular matrices $A,B \in \mathbb{R}^{n \times m}$, we want to find
a symmetric and positive semi-definite matrix $X^{*}\in \mathbb{R}^{n \times m} $ that solves the following optimization problem
\begin{equation}
  \min_{X\in \mathbb{R}^{n\times n}} \mathcal{F}(X) = \frac{1}{2}||XA - B||_F^{2} \quad s.t. \quad X\in\mathcal{S}_{+}(n), \label{problem}
\end{equation}
where $||Z||_F$ denotes the Frobenius norm of $Z\in \mathbb{R}^{r\times k}$ and $\mathcal{S}_{+}(n)$ represents the set of the symmetric and positive semi-definite $n$-by-$n$ matrices with real entries, i.e. $$\mathcal{S}_{+}(n) = \{X\in \mathbb{R}^{n\times n} : X^{\top} = X,\,\, v^{\top}Xv\geq 0,\,\forall v\in \mathbb{R}^{n}\}.$$

Problem (\ref{problem}) arises frequently in different applications such as: analysis of structures \cite{brock1968optimal, woodgate1996least}, signal processing \cite{suffridge1993approximation}, estimation of correlation matrices \cite{bagherpour2014efficient}, among other. It is well known that the feasible set $\mathcal{S}_{+}(n) $ of the problem (\ref{problem}), is a closed convex cone of dimensions $n\times (n+1)/2$, \cite{higham1988computing}. From this result and the fact that the objective function $\mathcal{F}$ is continuous, the existence of at least one global minimizer of problem (\ref{problem}) is guaranteed. Additionally, if $A$ is full rank then there exists a unique solution for (\ref{problem}), for more details see \cite{woodgate1996least}. In addition, since $\mathcal{F}$ is a convex function, this converts (\ref{problem}) in a convex minimization problem, which are relatively easy to solve. There are also some particular cases of the problem (\ref{problem}) that have an analytical formula for the solution, for example, when $A = I_{m}$ \cite{higham1988computing}, when $\verb"rank"(A) = 1 $ \cite{gillis2016semi}, or when $X$ is considered a diagonal matrix \cite{gillis2016semi}.\\

Although there are particular cases where problem (1) has a closed solution, in most cases, such as in real applications, a solution is only possible computationally, by an iterative method capable of dealing with non-stationary points, generating a sequence of feasible points that converge to a local minimizer. However, designing efficient algorithms that generate a feasible sequence of points is generally a difficult task, because it usually leads to use some projection operator, which requires computing spectral decompositions, which it is computationally expensive for large scale situations.\\

Due to the vast number of applications that problem (\ref{problem}) captures, many researchers are interested in studying this problem, both from a theoretical point of view, as well as from the design of new efficient approaches. In \cite{higham1988computing,liao2003least} the authors study a problem related to (\ref{problem}) theoretically, where they derive general analytical formulas for the solution considering some particular cases. Additionally, they provide sufficient and necessary conditions to guarantee the existence of the solution. On the other hand, the \emph{gradient projection method} was implemented by Nicolas Gillis et al. in \cite{gillis2016semi}. Specifically, Nicolas Gillis et.al. propose an algorithm called ``FGM'', which is an accelerated version of the classical gradient projection method, which uses the Nesterov \cite{nesterov2013introductory} acceleration technique. In essence, the FGM is an implementation of the algorithm that appears in \cite{nesterov2013introductory}. In \cite{gillis2016semi}, a method called ``AN-FGM'' is also proposed which is a semi-analytic approach that reduces the problem (\ref{problem}) to the case when $A$ is diagonal and then uses the FGM to address a more easy problem, this proposal looks quite efficient to deal with problems where $A$ is ill-conditioned. Another alternative to compute a numerical solution of problem (\ref{problem}) numerically, has been studied in \cite{andersson1997constrained}, where the authors propose an algorithm called ``Parallel tangents'' that is based on the gradient projection method that incorporates an over-relaxation step. One drawback of this parallel tangents method is that it does not guarantee optimal convergence. On the other hand, two algorithms that were designed to solve convex optimization problems over $\mathcal{S}_{+}(n) $, ``SDPT3'' \cite{toh1999sdpt3} and ``QSDP'' \cite{toh2008inexact} can be used to solve the problem (\ref{problem}).\\

In this work, we study the numerical behaviour of a spectral gradient projection method to address the positive semi-definite procrustes problem from a practical point of view. In particular, we adopt a gradient projection scheme with the non-monotone globalization technique proposed by Zhang and Hager in \cite{zhang2004nonmonotone}, in combination with the step size proposed by Barzilai and Borwein in \cite{barzilai1988two}. Subsequently, we present some computational experiments, in order to illustrate the effectiveness of the proposed method, solving the PSDP problem under many conditioning situations of $A$. Our main contribution is the implementation of an efficient gradient projection procedure using MATLAB, and the numerical comparison of the performance of such method against other existing algorithms of the state of the art.\\

The rest of this work is organized as follows. In the next section some important notations and tools for the good understanding of the article are introduced. In Section 3, the update formula of our proposed method is presented. Subsection 3.1, addresses the problem of selecting the step size of the proposed method and describes a non-monotone globalization technique to regulate such step size, and we culminate this subsection presenting the proposal. A derivation of the proposed method from the algorithm presented in \cite{francisco2017non} is discussed in subsection 3.2. In Section 4, several numerical experiments are carried out in order to demonstrate the effectiveness and efficiency of our procedure. Finally, the conclusions are presented in Section 5.

\section{Notations and Important Tools}
In this section, we present some fundamental concepts and tools for the well understanding of this work. Let's denote by $ \langle A,B \rangle := \sum_{i,j} A_{i,j} B_{i,j} = Tr[A^{\top}B]$ to the usual inner product on the matrix space $\mathbb{R}^{n\times m}$, here $Tr[A]$ denotes the trace of $A$. Given a differentiable function $\mathcal{F}: \mathbb{R}^{n\times m} \rightarrow \mathbb{R}$, the gradient of $\mathcal{F}$ respect to $X$ is denoted by $\nabla\mathcal{F}(X) = \left(\frac{\partial\mathcal{F}(X)}{\partial X_{i,j}}\right)$. The directional derivative of $\mathcal{F}$ at $X\in \mathbb{R}^{n \times m}$ in the direction $Z \in \mathbb{R}^{n \times m}$ is

\begin{equation}
\nabla\mathcal{F}(X)[Z] := \lim_{t\to\infty} \frac{\mathcal{F}(X+tZ) - \mathcal{F}(X)}{t} =  \langle \nabla\mathcal{F}(X),Z \rangle. \label{eq1}
\end{equation}

Another tool that we use is the projection operator over the feasible set $\mathcal{S}_{+}(n)$ which is defined below.
\begin{defi}
Let $X\in \mathbb{R}^{n\times n}$ be a real square matrix. The projection operator $\pi:\mathbb{R}^{n\times n} \rightarrow \mathcal{S}_{+}(n)$ over $\mathcal{S}_{+}(n)$ is defined by
\begin{equation}
\pi(X) = \textrm{arg}\min_{P\in \mathbb{R}^{n\times n}} ||P - X||_F, \quad s.t. \quad P\in\mathcal{S}_{+}(n). \label{eq2}
\end{equation}
\end{defi}

Note that the projection of any arbitrary matrix $X\in \mathbb{R}^{n \times n}$ is defined from an optimization problem, however, problem (\ref{eq2}) has a closed solution, this fact is established below.

\begin{prop}
Let $X\in \mathbb{R}^{n\times n}$ be a real square matrix. Then $\pi(X)$ is well-defined. Moreover, consider the symmetric part of $X$, that is, $X_{sym} = \frac{1}{2}(X^{\top} + X)$ and let $X_{sym} = V\Sigma V^{\top}$ be the spectral decomposition of $X_{sym}$, then $\pi(X) = V (\max(\Sigma,0)) V^{\top}$.
\end{prop}

\textbf{Proof.} The proof of this proposition appear in \cite{higham1988computing}.$\Box$

\section{A feasible Update Scheme}
Since $\mathcal{F}$ is smooth, a natural idea is to compute the next iterates as $Y(\tau) = X - \tau\nabla\mathcal{F}(X)$, where $X\in\mathcal{S}_{+}(n)$ is the previous iterate and $\tau>0$ represents the step size. The drawback of this approach is that the new point $Y(\tau)$ may not satisfies the constraints of problem (\ref{problem}). In order to overcome this issue, we consider the well-known \emph{projected gradient method} \cite{bertsekas1999nonlinear} which computes the new iterate $Z(\overline{\tau})$ as a point on the curve

\begin{equation}
Z(\tau) = \pi(X - \tau\nabla\mathcal{F}(X)). \label{eq3}
\end{equation}

Observe that equations (\ref{eq3}) guarantees that the new iterate preserves the feasibility. On the other hand, there are different techniques to select the step size $\tau$. The condition that is usually used for the \emph{gradient projection method} is known as ``Armijo's condition on the arc of projection'' \cite{bertsekas1999nonlinear}. This condition imposes to choose the step size $\tau_k$, at the $k$-th iteration, as the largest positive number $\tau$, such that verifies the following inequality
\begin{equation}
  \mathcal{F}(Z_k(\tau)) < \mathcal{F}(X_k)) + \sigma\nabla\mathcal{F}(X_k)[Z_k(\tau)-X_k], \label{Armijo}
\end{equation}
where $\sigma\in(0,1)$ and $Z_k(\tau) = \pi(X_k - \tau\nabla\mathcal{F}(X_k))$. The Armijo condition (\ref{Armijo}) is used in combination with a heuristic so-called \emph{backtracking} in order to find an appropriate step size that satisfies the condition (\ref{Armijo}), for more details about the backtracking strategy see \cite{bertsekas1999nonlinear,nocedal2006numerical}.\\

Note that if we ensure that the directional derivative $\nabla \mathcal{F}(X_k)[Z_k(\tau)-X_k]< 0$ for all $k$, then we obtain a sequence $\{X_k\}$ of points such that the corresponding sequence $\{\mathcal{F}(X_k)\}$ is monotonically decreasing. With the purpose of accelerating the convergence of the gradient projection scheme (\ref{eq3}), we adopt the non-monotone globalization technique proposed by Zhang and Hager in \cite{zhang2004nonmonotone} combined with the Barzilai and Borwein step sizes \cite{barzilai1988two} which usually accelerate the convergence of gradient-based methods. This strategy is described in the next section.

\subsection{The Barzilai-Borwein Step Sizes.}
In this subsection, we are focused on an non-monotone strategy for the step size selection as well as to present the proposed algorithm
in detail. \\

It is well-known that sometimes the Barzilai and Borwein step sizes \cite{barzilai1988two} can improve the performance of the gradient-based algorithms without increasing too much the computational cost of the procedure. Typically, this technique considers the steepest descent method, and proposes to choose any
of two step sizes, presente below, at the  $k$-th iteration,

\begin{eqnarray}
  \tau_{k}^{BB1} = \frac{||S_{k-1}||_F^{2}}{Tr[S_{k-1}^{\top}Y_{k-1}]}, &\textrm{or}& \tau_{k}^{BB2} = \frac{Tr[S_{k-1}^{\top}Y_{k-1}]}{||Y_{k-1}||_F^{2}}, \label{eq4}
\end{eqnarray}
where $S_{k-1} = X_{k} - X_{k-1}$ and $Y_{k-1} = \nabla \mathcal{F}(X_k) - \nabla \mathcal{F}(X_{k-1})$. Since the values $\tau_{k}^{BB1}$ and $\tau_{k}^{BB2}$ (BB-steps) could be negative, we used the absolute value of themselves to avoid negative step sizes that involve growth in the objective function. For more
details see \cite{barzilai1988two,raydan1997barzilai}. Since the BB-steps does not necessarily decrease the objective function values at each iteration, it can invalidate convergence. However, this issue can be overcome by using a globalization technique, which guarantees global convergence by regulating the step sizes in (\ref{eq4}), see \cite{dai2005projected,raydan1997barzilai}. Taking in mind this considerations, we adopt a non-monotone line
search method based on a strategy in \cite{zhang2004nonmonotone}, in our proposed algorithm. Specifically, the iterates
are recursively updated as $X_{k+1}: = Z_{k}(\tau_k) = \pi(X_k - \tau_k \nabla\mathcal{F}(X_k))$, where $\tau_k = \eta^{h} \tau_{k}^{BB1}$ or $\tau_k = \eta^{h} \tau_{k}^{BB2}$, where $h$ is the smallest integer number that verify the following condition

\begin{equation}
  \mathcal{F}(Z_{k}(\tau_k)) \leq C_k + \sigma \nabla \mathcal{F}(X_k)[Z_{k}(\tau_k)-X_k], \label{eq5}
\end{equation}
where $C_{k+1}$ is formed by the convex combination of $C_k$ and $\mathcal{F}(X_{k+1})$ given by $C_{k+1} = \frac{\mathcal{F}(X_{k+1}) + \gamma Q_kC_k}{Q_{k+1}}$, where $Q_{k+1} = \gamma Q_k + 1$, with $Q_0 = 1$. The proposed non-monotone gradient projection method to deal with the numerical solution of the problem (1) is summarized in Algorithm \ref{Alg1}.

\begin{algorithm}
\begin{algorithmic}[1]
\REQUIRE $X_{0}\in\mathcal{S}_{+}(n)$, $X_{-1} = X_0 + I_{n}$, $\tau>0$, $0<\tau_m\leq\tau_M$, $\sigma, \epsilon, \eta \in(0,1)$, $\gamma\in[0,1)$, $Q_0 = 1$, $C_0=\mathcal{F}(X_0)$, $k = 0$.
\WHILE{ $|| X_k - X_{k-1}) ||_F > \epsilon$ }
\WHILE { $\mathcal{F}(Z_{k}(\tau)) > C_k + \sigma\tau \mathcal{F}(X_k)[Z_{k}(\tau)-X_k]$}
\STATE $\tau = \eta\tau$,
\ENDWHILE
\STATE $X_{k+1} = Z_{k}(\tau)$, according to (\ref{eq3}).
\STATE Compute $Q_{k+1} = \gamma Q_k + 1$  and  $ C_{k+1} =(\gamma Q_kC_k + \mathcal{F}(X_{k+1}))/Q_{k+1}$.
\STATE Take $\tau = |\alpha_{k+1}^{BB1}|$ or well $\tau = |\alpha_{k+1}^{BB2}|$, according to (\ref{eq4}).
\STATE $\tau = \max(\min(\tau,\tau_M),\tau_m)$.
\STATE $k = k+1$.
\ENDWHILE
\STATE $X^{*} = X_k$.
\end{algorithmic}
\caption{OptPSDP}\label{Alg1}
\end{algorithm}

\begin{obs}
Note that if we select $\gamma = 0$ in the previous algorithm, then Algorithm \ref{Alg1} is reduced to the classical gradient projection method. Observe also that the Algorithm \ref{Alg1} can be used to minimize any objective smooth function over the matrix set $\mathcal{S}_{+}(n) $, however, the interest of this work is focused on the particular problem (\ref{problem}).
\end{obs}

Note that the step 5, in Algorithm \ref{Alg1}, is the step that require the major computational effort, because it needs to use the operator
projection defined in (\ref{eq2}), which requires compute a spectral decomposition, which it is computationally inefficient. In order to avoid
the calculation of such spectral decomposition, in each step, we propose the following idea: first note that if the symmetric part of $Y_k = X_k - \tau_k \nabla\mathcal{F}(X_k)$ is positive definite then this matrix coincides with its projection over $\mathcal{S}_{+}(n)$. Thus, we propose to use the Cholesky's factorization to make the Algorithm \ref{Alg1} more efficient. Specifically, in the step 5, we try to compute the Cholesky factorization of $\frac{Y_k + Y_k^{\top}}{2}$, if no error is generated, then $X_{k+1}$ is updated by $X_{k+1} = \frac{Y_k+Y_k^{\top}}{2}$, otherwise $X_{k+1} = Z_{k}(\tau_k)$ is updated using the projection operator. In Section \ref{experiments}, we demostraste numerically the efficiency of this strategy on some numerical tests.

\subsection{Another Point of View of Algorithm \ref{Alg1}.}
In this section we derive Algorithm \ref{Alg1} from an algorithm proposed by Francisco et al. in \cite{francisco2017non} recently. In addition, we establish a convergence result related to our Algorithm \ref{Alg1}.\\

In \cite{francisco2017non} the authors propose a globally convergent non-monotonous algorithm to numerically solve the following optimization problem,
\begin{equation}
  \min f(x) \quad s.t. \quad x\in \Omega, \label{problem2}
\end{equation}
where $\Omega$ is a closed subset of $\mathbb{R}^{n}$ and $f: \mathbb{R}^{n}\rightarrow \mathbb{R}$ is a continuously differentiable function on $\hat{\Omega}$ such that $\Omega \subset \hat{\Omega}$. The proposed algorithm by Francisco et.al. builds a sequence of iterates as follows: given the current point $x_k \in\Omega$, $\rho_k> 0$ a positive scalar and two symmetric matrices $A_k, B_k$, with $A_k$ definite positive, then the next trial point $x_{k+1}$ is computed as the argument that minimize the quadratic model
\begin{equation}
\min_{x\in\Omega}  Q_{k}(x) = \langle\nabla f(x_k),x - x_k\rangle + \frac{1}{2} (x - x_k)^{\top}(B_k + \rho_k A_k)(x - x_k), \label{model}
\end{equation}
where $\rho_k$ works as a regularization parameter. This method is based on the ideas of the trust region methods \cite{nocedal2006numerical} and the well-known method of Levenberg-Marquardt \cite{nocedal2006numerical}. The authors in \cite{francisco2017non}, combine these ideas with the non-monotone technique proposed by Zhang and Hager \cite{zhang2004nonmonotone}, and thus obtain a very general method to solve the non-linear optimization problem (\ref{problem2}). \\

The rest of this subsection is dedicated to demonstrate that the Algorithm \ref{Alg1} can be seen as a particular case of the algorithm proposed in \cite{francisco2017non}. For this end, it is sufficient to demonstrate that the update formula of our proposal (\ref{eq3}) is equivalent to solve a quadratic model on
$\mathcal{S}_{+}(n)$, due to the non-monotone strategy to choose the step size is the same for the two algorithms.

\begin{prop}\label{prop2}
   Let $X_k\in\mathcal{S}_{+}(n)$ be the point generated by Algorithm \ref{Alg1} at the $k$-th iteration. If $\tau>0$ then $Z_k(\tau) = \pi(X_k - \tau\nabla\mathcal{F}(X_k))$ is the minimum of the following quadratic model,
  \begin{equation}
  \mathcal{Q}_{k}(X) = Tr[\nabla \mathcal{F}(X_k)^{\top}(X - X_k)] + \frac{1}{2\tau}||X - X_k||_F^{2}, \label{model2}
  \end{equation}
  over the set $\mathcal{S}_{+}(n)$.
\end{prop}

\textbf{Proof.} Since $Z_k(\tau) = \pi(X_k - \tau\nabla\mathcal{F}(X_k))$ then $Z_k(\tau)$ is a solution of
  \begin{equation}
    \min \mathcal{J}(X) = \frac{1}{2}||X - ( X_k - \tau\nabla\mathcal{F}(X_k) )||_F^{2}, \quad s.t. \quad P\in\mathcal{S}_{+}(n), \label{J_fun}
  \end{equation}

   From the definition of $J(X)$ and using trace properties we have
   \begin{eqnarray}
     \mathcal{J}(X) & = & \frac{1}{2} Tr[( X - (X_k - \tau\nabla\mathcal{F}(X_k)))^{\top}( X - (X_k - \tau\nabla\mathcal{F}(X_k)))], \nonumber \\
          & = & \frac{1}{2} Tr[X^{\top}X - 2X^{\top}X_k + 2\tau X^{\top}\nabla\mathcal{F}(X_k) + X_k^{\top}X_k - 2\tau X_k^{\top}\nabla\mathcal{F}(X_k)  \nonumber \\
          &   & + \tau^{2}\nabla\mathcal{F}(X_k)^{\top}\nabla\mathcal{F}(X_k)].
   \end{eqnarray}

  Now, since $\tau^{2}\nabla\mathcal{F}(X_k)^{\top}\nabla\mathcal{F}(X_k)$ is constant, then minimize $\mathcal{J}(\cdot)$ is equivalent to minimize the function $\hat{\mathcal{J}}(\cdot)$ given by
  $$\hat{\mathcal{J}}(X) =  \frac{1}{2} Tr[X^{\top}X - 2X^{\top}X_k + 2\tau X^{\top}\nabla\mathcal{F}(X_k) + X_k^{\top}X_k - 2\tau X_k^{\top}\nabla\mathcal{F}(X_k)], $$

  Rewriting this last result we arrive at
  \begin{equation}
    \hat{\mathcal{J}}(X) = \tau Tr[X^{\top}\nabla\mathcal{F}(X_k) - X_k^{\top}\nabla\mathcal{F}(X_k)] + \frac{1}{2}Tr[X^{\top}X - 2X^{\top}X_k + X_k^{\top}X_k],
  \end{equation}
  or equivalently
  \begin{equation}
    \hat{\mathcal{J}}(X) = \tau\left(  Tr[\nabla\mathcal{F}(X_k^{\top})(X-X_k)] + \frac{1}{2\tau}||X-X_k||_F^{2}  \right).
  \end{equation}
  Then, since $\tau$ is constant for the optimization process over $\mathcal{S}_{+}(n)$, we have that minimize $\hat{\mathcal{J}}(\cdot)$ over $\mathcal{S}_{+}(n)$, is equivalent to minimize the quadratic function $\mathcal{Q}_{k}(X)$ defined in (\ref{model2}) over the set $\mathcal{S}_{+}(n)$, which completes the proof.$\Box$\\

Note that Proposition \ref{prop2} shows that the Algorithm \ref{Alg1} is a particular case of the algorithm proposed by Francisco et al. \cite{francisco2017non}, obtained  taking $\rho_k = 1$, $B_k$ to the null matrix and $A_k = \frac{1}{\tau_k}I_n$ at each iteration. This result implies that the Algorithm \ref{Alg1} is globally convergent, which it is establishes in Theorem \ref{teorema1}.

\begin{teo}\label{teorema1}
Let $\{X_k\}$ be a sequence generated by Algorithm \ref{Alg1}. Assume that $\gamma<1$, then every accumulation points of $\{X_k\}$ is a stationary point of the problem (\ref{problem}).
\end{teo}

\section{Numerical Experiments}
\label{experiments}
In this section, we illustrate the effectiveness and efficiency of the proposed algorithm (Algorithm 1: OptPSDP) on several positive semi-definite procrustes problems generated synthetically. An implementation in Matlab of OptPSDP is available at: \textsf{http://www.mathworks.com/matlabcentral
/fileexchange/64597-spectral-projected-gradient
-method-for-the-positive-semi-definite-procrustes-problem}.\\

All computational experiments were carried out using Matlab 7.0 in an intel (R) CORE (TM) i7-4770 processor, 3.40 GHz CPU with 500 Gb of HD and 16 Gb of Ram. In all experiments the following values are used for the OptPSDP algorithm: $\sigma = $ 1e-4, $\tau_0 = $ 1e-3, $\tau_{\min} = $ 1e-20, $ \tau_{\max} = $ 1e20, $\epsilon = $ 1 e-5, $\gamma = 0.85$ and $\eta = 0.2$. For the other methods we used the default parameters of each algorithm, except for the tolerance fixed to $\epsilon = $ 1e-5. As a maximum number of iterations $N = 10000$ was selected for all algorithms.\\

\begin{table}
\centering
\caption{Numerical results for well conditioned PSDP (problem = 1).}
{\begin{tabular}{lcccccc}
 \hline
\textbf{Method}	            &	\textbf{Nitr}	&	\textbf{Nfe}	&	\textbf{Time}	&	\textbf{XErr}	&	\textbf{Fval}	&	\textbf{Global Error}\\
\hline
&\multicolumn{6}{c}{E1: n = 100, m = 70, problem = 1, $\gamma$ = 0.8 }\\
\cline{2-7}
\textbf{Grad}	        &	746	    &	747	    &	0.808	&	9.96e-6	&	1.06e-6	    &	17.3	\\
\textbf{FGM}	        &	1756	&	1757	&	1.889	&	9.99e-6	&	7.43e-25	&	61.3	\\
\textbf{ParTan}	        &	74 	    &	75	    &	\textbf{0.107}	&	2.91e-6	&	1.02e-8	    &	17.4	\\
\textbf{OptPSDP}	    &	102	    &	103	    &	0.123	&	8.21e-6	&	5.86e-7	    &	17.3	\\
\hline
& \multicolumn{6}{c}{E2: n = 150, m = 100, problem = 1, $\gamma$ = 0.85}	\\
\cline{2-7}
\textbf{Grad}	        &	1356	&	1357	&	3.166	&	9.98e-6	&	2.29e-6	    &	27.7	\\
\textbf{FGM}	        &	2091	&	2092	&	4.856	&	9.99e-6	&	2.37e-24	&	101.2	\\
\textbf{ParTan}	        &	101	    &	102	    &	\textbf{0.316}	&	2.38e-6	&	1.17e-8	    &  	27.8	\\
\textbf{OptPSDP}	    &	153	    &	154	    &	0.393	&	8.68e-6	&	1.42e-6  	&	27.6	\\
\hline
& \multicolumn{6}{c}{E3: n = 1000, m = 100, problem = 1, $\gamma$ = 0.85}\\
\cline{2-7}
\textbf{Grad}	        &	8	&	9	&	1.656	&	4.91e-6	&	1.60e-11	&	4.76e-7	\\
\textbf{FGM}	        &	8	&	9	&	1.588	&	2.68e-6	&	2.22e-12	&	1.31e-7	\\
\textbf{ParTan}	        &	7	&	8	&	2.287	&	3.83e-6	&	9.28e-12	&	3.70e-7	\\
\textbf{OptPSDP}	    &	7	&	8	&	\textbf{1.438}	&	5.86e-6	&	5.79e-12	&	2.72e-7	\\
\hline
&\multicolumn{6}{c}{E4: n = 1500, m = 1500, problem = 1, $\gamma$ = 0.85}\\
\cline{2-7}
\textbf{Grad}	        &	8	&	9	&	5.471	&	4.75e-6	&	1.65e-11	&	4.52e-7	\\
\textbf{FGM}	        &	8	&	9	&	5.261	&	1.15e-6	&	5.64e-13	&	7.16e-8	\\
\textbf{ParTan}	        &	7	&	8	&	7.261	&	5.15e-6	&	1.64e-11	&	4.90e-7	\\
\textbf{OptPSDP}	    &	8	&	9	&	\textbf{4.752}	&	2.72e-6	&	2.41e-12	&	1.24e-7	\\
\hline
\end{tabular}}\label{tab:1}
\end{table}

In the rest of this section, we denote by ``Nitr'' the average number of iterations, ``Nfe'' the average number of functions evaluations, ``Time'' the average execution time in seconds, ``Fval'' the average value of the evaluation of the objective function at point $\hat{X}$ which denotes the optimum estimated by each algorithm, ``Error'' the average global error, that is, $||X^{*} - \hat{X}||_F$, where $X^{*}$ denotes the global optimum each PSDP problem, and finally we denote by ``XErr'', the average error $||\hat{X}-X_{k}||_F$, and $X_ {k}$ penultimate point generated by each algorithm. In addition, we denote by \textbf{Grad} to the classical gradient projection method proposed in \cite{gillis2016semi}, \textbf{FGM} denotes the accelerated gradient projection method proposed in \cite{gillis2016semi}, \textbf{ParTan} denotes parallel tangent method introduced in \cite{andersson1997constrained} and \textbf{OptPSDP} denotes our proposal.\\

For the numerical experiments, we consider problem (\ref{problem}) where the matrix $A\in \mathbb{R}^{n \times m}$ is build as $A = P\Lambda Q^{\top}$, where $P\in \mathbb{R}^{n \times n}$ and $Q\in \mathbb{R}^{m \times m}$ are orthogonal matrices randomly generated and $\Lambda \in \mathbb{R}^{n \times m}$ is a diagonal matrix defined as we explain below. The starting point $X_0$  was generated as $ X_0 = \pi(\bar{X}_0)$, where $\bar{X}_0$ was randomly generated. In order to monitoring the behavior of the algorithms, the optimal solution is generated by $X^{*} = \pi(\tilde{X})$ where $\tilde{X} \in \mathbb{R}^{n \times n}$ was randomly generated. Then, the matrix $B \in \mathbb{R}^{n \times m}$ was taken as $B = XA$, in this way,  $X^{*}$ is a global optimum of the problem (\ref{problem}) with optimal value zero, i.e. $\mathcal{F}(X^{*}) = 0$. All random values were generated following a standard normal distribution using the \verb"randn" function of Matlab.\\

\begin{table}
\centering
\caption{Numerical results for ill conditioned PSDP (problem = 2).}
{\begin{tabular}{lcccccc}
 \hline
\textbf{Method}	            &	\textbf{Nitr}	&	\textbf{Nfe}	&	\textbf{Tiempo}	&	\textbf{XErr}	&	\textbf{Fval}	&	\textbf{Global Error}\\
\hline
&\multicolumn{6}{c}{E5: n = 30, m = 10, problem = 2, $\gamma$ = 0.85}\\
\cline{2-7}
\textbf{Grad}	        &	4392	&	4393	&	0.81	&	1.15e-6	&	6.2e-5	&	9.7	\\
\textbf{FGM}	        &	1289	&	1290	&	0.251	&	9.96e-6	&	3.55e-9	&	16.5	\\
\textbf{ParTan}	        &	187	    &	188	    &	\textbf{0.041}	&	5.04e-6	&	3.49e-6	&	9.8	\\
\textbf{OptPSDP}	    &	392	    &	395	    &	0.086	&	8.5e-6	&	1.22e-5	&	9.7	\\
\hline
&\multicolumn{6}{c}{E6: n = 100, m = 50, problem = 2, $\gamma$ = 0.55}\\
\cline{2-7}
\textbf{Grad}	        &	10000	&	10001	&	12.123	    &	3.4e-5	&	1.41e-2	&	25.3	\\
\textbf{FGM}	        &	5181	&	5182	&	6.4699	    &	9.85e-6	&	5.1e-6	&	54.9	\\
\textbf{ParTan}	        &	1073	&	1074	&	\textbf{1.553}	    &	8.32e-6	&	1.27e-4	&	25.6	\\
\textbf{OptPSDP}	    &	1221	&	1231	&	1.57	    &   9.68e-6	&	2.4e-3	&	25.3	\\
\hline
&\multicolumn{6}{c}{E7: n = 60, m = 60, problem = 2, $\gamma$ = 0.85}\\
\cline{2-7}
\textbf{Grad}	        &	6149	&	6150	&	2.731	&	1.16e-6	&	6.97e-4	&	1.93e-2	\\
\textbf{FGM}	        &	442	    &	443	    &	0.236	&	9.86e-6	&	2.04e-7	&	3.77e-4	\\
\textbf{ParTan}	        &	291	    &	292	    &	0.196	&	4.63e-6	&	2.51e-6	&	1.1e-3	\\
\textbf{OptPSDP}	    &	318	    &	322	    &	\textbf{0.184}	&	8.6e-6	&	1.24e-4	&	9.3e-3	\\
\hline
&\multicolumn{6}{c}{E8: n = 120, m = 120, problem = 2, $\gamma$ = 0.55}\\
\cline{2-7}
\textbf{Grad}	        &	9962	&	9963	&	14.679	&	4.55e-6	&	1.46e-1	&	2.58e-1	\\
\textbf{FGM}	        &	784	    &	785	    &	1.394	&	9.93e-6	&	8.4e-7	&	6.82e-4	\\
\textbf{ParTan}	        &	498	    &	499	    &	\textbf{1.14}	&	5.1e-6	&	1.68e-5	&	2.7e-3	\\
\textbf{OptPSDP}	    &	728	    &	743	    &	1.387	&	9.08e-6	&	1.7e-3	&	3.11e-2	\\
\hline
\end{tabular}}\label{tab:2}
\end{table}

In addition, we consider the following three distributions of the entries of $\Lambda$,\\

\textbf{Problema 1}: The $\Lambda$ diagonal entries are generated by a truncated normal distribution in the interval [10,12].\\

\textbf{Problema 2}: The diagonal of $\Lambda$ is given by $\lambda_{ii} = i + 2r_{i}$, where $r_i$ is a randomly generated from the uniform distribution in the interval [0,1].\\

\textbf{Problema 3}: Each element of the diagonal matrix $\Lambda$ is generated as $\lambda_{ii} = 1 + \frac{99 (i-1)}{m+1} + 2r_{i}$, with $r_i$ is a randomly generated from the uniform distribution in the interval [0,1].\\

Observe that if the $\Lambda$ is generated following the structure of \textbf{Problema 1} then $A$ is a well-conditioned matrix, while it is generated by the diagonal
structures describe in \textbf{Problema 2} and \textbf{Problema 3} then $ A $ is a ill-conditioned matrix. In order to study the numerical behavior and performance of all methods,  we  consider several size of problems PSDP and different conditions number of $A$. In all tables, we present the averages of the comparing values obtained by each algorithm in a total of 50 independent instances.\\

In the first experiment, we study the efficiency of the proposed method on well-conditioned PSDP problems. Table \ref{tab:1} summarizes the numerical results of this comparison. From Table \ref{tab:1} we observe that the methods that the methods that converge faster are \textbf{ParTan} and \textbf{OptPSDP}. In addition, it's seen that if $A$ is rectangular then the most efficient method in terms of CPU-time is \textbf{ParTan}. However, clearly we note that our proposal is more efficient for problems where $A$ is square. According to the error \emph{XErr}, all algorithms reach an order less than 1e-5 and additionally, we can see that the value \emph{Fval} is close to zero for all algorithms.\\

\begin{table}
\centering
\caption{Numerical results for ill conditioned PSDP (problem = 3).}
{\begin{tabular}{lcccccc}
 \hline
\textbf{Method}	            &	\textbf{Nitr}	&	\textbf{Nfe}	&	\textbf{Time}	&	\textbf{XErr}	&	\textbf{Fval}	&	\textbf{Global Error}\\
\hline
&\multicolumn{6}{c}{E9: n = 50, m = 10, problem = 3, $\gamma$ = 0.55}\\
\cline{2-7}
\textbf{Grad}	        &	10000	&	10001	&	3.539	&	6.17e-5	&	1.45e-1	&	19.51	\\
\textbf{FGM}	        &	6233	&	6234	&	2.287	&	1.08e-6	&	1.11e-5	&	26.29	\\
\textbf{ParTan}	        &	514	    &	515	    &	\textbf{0.21}	&	7.06e-6	&	3.89e-4	&	19.53	\\
\textbf{OptPSDP}	    &	1867	&	1885	&	0.734	&	9.62e-6	&	4.1e-3	&	19.68	\\
\hline
&\multicolumn{6}{c}{E10: n = 100, m = 10, problem = 3, $\gamma$ = 0.55}\\
\cline{2-7}
\textbf{Grad}	        &	10000	&	10001	&	10.299	&	1.2e-4	&	5.73e-1	&	43.44	\\
\textbf{FGM}	        &	7108	&	7109	&	7.951	&	1.39e-6	&	6.86e-6	&	48.15	\\
\textbf{ParTan}	        &	626	    &	627	    &	\textbf{0.791}	&	5.21e-6	&	1.19e-4	&	43.3	\\
\textbf{OptPSDP}	    &	3817	&	3840	&	4.206	&	1.04e-6	&	9.7e-3	&	43.58	\\
\hline
&\multicolumn{6}{c}{E11: n = 100, m = 100, problem = 3, $\gamma$ = 0.55}\\
\cline{2-7}
\textbf{Grad}	        &	9250	&	9251	&	9.625	&	3.14e-6	&	4.89e-2	&	1.52e-1	\\
\textbf{FGM}	        &	686	    &	687	    &	0.86	&	9.9e-6	&	5.56e-7	&	5.91e-4	\\
\textbf{ParTan}	        &	437	    &	438	    &	\textbf{0.7}	    &	4.72e-6	&	8.88e-6	&	2.2e-3	\\
\textbf{OptPSDP}	    &	606	    &	611	    &	0.801	&	8.94e-6	&	9.07e-4	&	2.42e-2	\\
\hline
&\multicolumn{6}{c}{E12: n = 150, m = 150, problem = 3, $\gamma$ = 0.55}\\
\cline{2-7}
\textbf{Grad}	        &	9660	&	9661	&	21.894	&	4.17e-6	&	7.5e-2	&	2.04e-1	\\
\textbf{FGM}	        &	729	    &	730	    &	1.986	&	9.91e-6	&	5.63e-7	&	6.2e-4	\\
\textbf{ParTan}	        &	486	    &	487	    &	\textbf{1.764}	&	2.9e-6	&	3.68e-6	&	1.3e-3	\\
\textbf{OptPSDP}	    &	649	    &	665	    &	1.896	&	8.87e-6	&	1e-3	&	2.66e-2	\\
\hline
\end{tabular}}\label{tab:3}
\end{table}

In tables \ref{tab:2} and \ref{tab:3} we present the results obtained by the four procedures solving ill-conditioned PSDP. These tables clearly show that \textbf{Grad} algorithm is the method that obtain the worst results, because sometimes run the maximum number of iterations allowed and it is the slowest in terms CPU-time. On the other hand, we observe that the \textbf{FGM}, \textbf{ParTan} and \textbf{OptPSDP} methods show similar performance both in the number of iterations, and in execution time when $ m = n $. However, when $A$ is a rectangular matrix, the more efficient method is \textbf{ParTan}. In spite of this, all the methods reach convergence, since all obtain  small values of \emph{XErr}.\\

For the fourth experiment group, the PSDP problems were constructed with randomly generated synthetic data as explained at the beginning of this section, however, the optimum $X^{*}$ matrix was built as follows, first a matrix $M\in \mathbb{R}^{n \times n}$ is randomly generated with entries following a standard normal distribution, afterwards $V$ is obtained as the orthogonal matrix of the QR factorization of $M$, from this matrix, we set $X^{*} = V^{\top}\Sigma V $, where $\Sigma \in \mathbb{R}^{n \times n}$ is a diagonal matrix whose diagonal elements were generated by $\Sigma(1,1) = \Sigma(2,2) = 0$ and $\Sigma(i,i) = \verb"rand"$ for all $i\in\{3,4,\ldots,n\}$ using Matlab notation. Thus, the optimal solution of the PSDP generated is a symmetric and positive semi-definite matrix with only two eigenvalues equal to zero and $n-2$ strictly positive eigenvalues.\\

\begin{table}
\centering

\caption{Numerical results for several PSDP (problem = 1,2,3).}
{\begin{tabular}{lcccccc}
 \hline
\textbf{Method}	            &	\textbf{Nitr}	&	\textbf{Nfe}	&	\textbf{Time}	&	\textbf{XErr}	&	\textbf{Fval}	&	\textbf{Global Error}\\
\hline
&\multicolumn{6}{c}{E13: n = 10, m = 70, problem = 1, $\gamma$ = 0.55}\\
\cline{2-7}
\textbf{Grad}	        &	25	&	26	&	0.031	&	7.61e-6	&	5.02e-9	    &	18.85	\\
\textbf{FGM}	        &	788	&	789	&	0.956	&	9.98e-6	&	5.73e-26	&	22.27	\\
\textbf{ParTan}	        &	12	&	13	&	0.019	&	4.38e-6	&	1.7e-10	    &	18.91	\\
\textbf{OptPSDP}	    &	11	&	12	&	\textbf{0.009}	&	6.00e-6	&	2e-9	    &	18.82	\\
\hline
&\multicolumn{6}{c}{E14: n = 1000, m = 1000, problem = 1, $\gamma$ = 0.85}\\
\cline{2-7}
\textbf{Grad}	        &	8	&	9	&	1.868	&	4.57e-6	&	1.56e-11	&	4.73e-7	\\
\textbf{FGM}	        &	8	&	9	&	1.825	&	1.86e-6	&	1.51e-12	&	9.8e-8	\\
\textbf{ParTan}	        &	7	&	8	&	2.482	&	3.76e-6	&	8.44e-12	&	3.53e-7	\\
\textbf{OptPSDP}	    &	8	&	9	&	\textbf{0.935}	&	2.88e-6	&	3.09e-12	&	1.45e-7	\\
\hline
&\multicolumn{6}{c}{E15: n = 60, m = 30, problem = 2, $\gamma$ = 0.85}\\
\cline{2-7}
\textbf{Grad}	        &	3737	&	3738	&	1.892	&	9.99e-6	&	2.97e-5	&	15.9	\\
\textbf{FGM}	        &	1876	&	1877	&	0.9692	&	9.22e-6	&	2.39e-6	&	16.4	\\
\textbf{ParTan}	        &	185	    &	186	    &	0.113	&	3.43e-6	&	3.6e-7	&	15.9	\\
\textbf{OptPSDP}	    &	178	    &	183	    &	\textbf{0.063}	&	7.49e-6	&	8.42e-6	&	15.9	\\
\hline
&\multicolumn{6}{c}{E16: n = 100, m = 100, problem = 2, $\gamma$ = 0.55}\\
\cline{2-7}
\textbf{Grad}	        &	9437	&	9438	&	11.787	&	3.18e-6	&	5.76e-2	&	1.57e-1	\\
\textbf{FGM}	        &	642	    &	643	    &	0.814	&	9.91e-6	&	5.84e-7	&	5.78e-4	\\
\textbf{ParTan}	        &	378	    &	379	    &	0.601	&	4.2e-6	&	8.99e-6	&	1.6e-3	\\
\textbf{OptPSDP}	    &	511	    &	536	    &	\textbf{0.47}	&	8.27e-6	&	8.05e-4	&	2.13e-2	\\
\hline
&\multicolumn{6}{c}{E17: n = 60, m = 30, problem = 3, $\gamma$ = 0.55}\\
\cline{2-7}
\textbf{Grad}	        &	10000	&	10001	&	5.073	&	6.14e-5	&	1.77e-1	&	15.8	\\
\textbf{FGM}	        &	4179	&	4180	&	2.168	&	9.49e-6	&	6.1e-5	&	16.3	\\
\textbf{ParTan}	        &	596	    &	597	    &	\textbf{0.365}	&	6.54e-6	&	3.87e-5	&	15.9	\\
\textbf{OptPSDP}	    &	1366	&	1390	&	0.429	&	8.69e-6	&	1.1e-3	&	17.8	\\
\hline
&\multicolumn{6}{c}{E18: n = 120, m = 120, problem = 3, $\gamma$ = 0.55}\\
\cline{2-7}
\textbf{Grad}	        &	9281	&	9282	&	16.522	&	4.21e-6	&	8.85e-2	&	2.13e-1	\\
\textbf{FGM}	        &	672	    &	673	    &	1.218	&	9.91e-6	&	5.64e-7	&	6.01e-4	\\
\textbf{ParTan}	        &	393	    &	394	    &	0.897	&	3.07e-6	&	6.98e-6	&	1.4e-3	\\
\textbf{OptPSDP}	    &	557	    &	582	    &	\textbf{0.686}	&	8.29e-6	&	8.56e-4	&	2.32e-2	\\
\hline
\end{tabular}}\label{tab:4}
\end{table}

The numerical results corresponding to the third experiment are shown in Table \ref{tab:4}. This table shows that the \textbf{ParTan} algorithm obtained the best performance in terms of the number of iterations, in almost all experiments. In addition, we observe that the our \textbf{OptPSDP} is the most efficient procedure in terms of CPU-time in both well-conditioned and ill-conditioned problems. From all the experiments performed, we concluded that the our proposal is a competitive alternative to solve the problem \ref{problem} under different situations of conditioning and scale of $A$.

\section{Conclusions}
The problem (\ref{problem}) has a wide range of applications in the fields of structure analysis, physical problems, signal processing, estimation of correlation matrices, among others. To address this problem, we design and implement an efficient and globally convergent algorithm that preserves feasibility in each iteration. Our proposal is based on the gradient projection method and we incorporate a non-monotone strategy in combination with the Barzilai and Borwein step sizes in order to accelerate the convergence. The bottleneck of the proposed algorithm is the computation of the projection operator, which is computationally inefficient. In order to improve the efficiency of our algorithm, we present a strategy based on Cholesky factorization to reduce the number of projections. This technique can be a good alternative to deal with large-scale problems. Some theoretical results were presented. Finally, from the numerical experiments we note that the performance of the resulting algorithm is quite competitive with some of the state of the art methods.

\paragraph{Acknowledgements:} 
This research was supported in part by Conacyt, Mexico (PhD. studies scholarship).

\bibliographystyle{plain}
\bibliography{journal0}

\end{document}